\documentclass[times, 11pt]{article}
\usepackage{epsfig,amsfonts,amssymb, amsmath}

\newtheorem{theorem}{Theorem}[section]
\newtheorem{lemma}[theorem]{Lemma}
\newtheorem{corollary}[theorem]{Corollary}

\newtheorem{proposition}[theorem]{Proposition}
\newtheorem{definition}[theorem]{Definition}

\newtheorem{Def}[theorem]{Definition}

\newcommand{\proof}{
{\bf Proof: }}

\newcommand{\eat}[1]{}

\newcommand{\qed}{\nobreak \ifvmode \relax \else
          \ifdim\lastskip<1.5em \hskip-\lastskip
          \hskip1.5em plus0em minus0.5em \fi \nobreak
          \vrule height0.75em width0.5em depth0.25em\fi \newline}

\newcommand{\Sum}{
\displaystyle\sum}

\newcommand{\Prod}{
\displaystyle\prod}

\newcommand{\mc}[1]{{\ensuremath \mathcal{#1}}}
\newcommand{\zo}{\{0,1\}}

\newcommand{\etal}{{\em et al.}}
\newcommand{\bp}{\{-1,+1\}}
\newcommand \PP {{\rm PP}}
\newcommand \Par {{\rm PAR}}
\newcommand \IP {{\rm IP}}
\newcommand \AND {{\rm AND}}
\newcommand \Or {{\rm OR}}
\newcommand \Th {{\rm THR}}
\newcommand \Thr {{\rm THR {\ensuremath \circ} AND}}
\newcommand \AC {{\ensuremath {\mathrm AC^0}}} 

\newcommand \R {\mathbb{R}}
\newcommand \Z {\mathbb{Z}}

\newcommand \Xn {{\ensuremath X_1,\ldots, X_n}}
\newcommand \Yn {{\ensuremath Y_1,\ldots, Y_n}}
\newcommand \an {{\ensuremath a_1,\ldots,a_n}}
\newcommand \Xno {{\ensuremath X_1,\ldots,X_{n-1}}}
\newcommand \Yno {{\ensuremath Y_1,\ldots,Y_{n-1}}}

\newcommand{\ol}[1]{\ensuremath \overline{#1}}
\newcommand{\sprb}{{\ensuremath \spr_{\mc{B}}}}

\DeclareMathOperator{\spr}{sp}
\DeclareMathOperator{\sgn}{sgn}

\begin{document}

\title{Polynomials that Sign Represent Parity\\ and\\ Descartes' Rule of Signs}
\author{
        Saugata Basu\thanks{Supported in part by an NSF Career Award 0133597 and a Sloan Foundation        
Fellowship}\\
        School of Mathematics\\
        Georgia Tech\\
        {\tt saugata@math.gatech.edu}
\and    Nayantara Bhatnagar\\
        College of Computing\\
        Georgia Tech\\
        {\tt nand@cc.gatech.edu}
 \and Parikshit Gopalan\\
        College of Computing\\
        Georgia Tech\\
        {\tt parik@cc.gatech.edu}
 \and Richard J. Lipton \thanks{Also with Telcordia. Supported in part by NSF CCR-0002299}\\
        College of Computing\\
        Georgia Tech\\
        {\tt rjl@cc.gatech.edu}
}

\maketitle
\thispagestyle{empty}

\begin{abstract}

A real polynomial $P(X_1,\ldots, X_n)$ sign represents $f: A^n \rightarrow \zo$ if
for every $(a_1, \ldots, a_n) \in A^n$, the sign of
$P(a_1,\ldots,a_n)$ equals $(-1)^{f(a_1,\ldots,a_n)}$. Such sign
  representations are well-studied in computer science and have
 applications to computational complexity and computational learning
 theory. The work in this area aims to determine the minimum degree and
 sparsity possible for a polynomial that sign represents
 a function $f$. While the degree of such polynomials is
 relatively well-understood, far less is known about their
 sparsity. Known bounds apply only to the cases where $A =\zo$ or $A =
 \bp$.

In this work, we present a systematic study of tradeoffs between
degree and sparsity of sign representations through the lens of the
parity function. We attempt to prove
bounds that hold for any choice of set $A$. 
We show that sign representing parity over $\{0,\ldots,
  m-1\}^n$ with the degree in each variable at most $m-1$ requires
  sparsity at least $m^n$. 
We show that a tradeoff exists between sparsity and degree, by
exhibiting a sign representation that has higher degree but lower
  sparsity. We show a lower bound of $n(m -2) + 1$ on the sparsity
  of polynomials of any degree representing parity over $\{0,
  \ldots, m-1\}^n$. We prove exact bounds on the sparsity of such polynomials for
  any two element subset $A$. The main tool used is Descartes' Rule of Signs, a
classical result in algebra,  relating the sparsity of a polynomial to
its  number of real roots.

As an application, we use bounds on sparsity to derive circuit lower
  bounds for depth-two AND-OR-NOT circuits with a Threshold
  Gate at the top. We use this
  to give a simple proof that such circuits need size $1.5^n$ 
  to compute parity, which improves the previous bound of $\frac{4}{3}^{n/2}$ \cite{Goldmann}. We
  show a tight lower bound of $2^n$ for the inner product function over $\{0,1\}^n \times \{0, 1\}^n$.

\end{abstract}
\newpage

\section{Introduction}

Let $A$ be a subset of $\Z$ and let $f: A^n \rightarrow \{0,1\}$
be a function on $A^n$. 

\begin{Def}
\label{def:exact}
A polynomial $P(X_1,\ldots, X_n) \in \R[X_1,\cdots,X_n]$ exactly represents $f$ over $A^n$ if for every
$(a_1, \ldots, a_n) \in A^n$, $P(a_1,\ldots,a_n) = f(a_1, \ldots, a_n)$.
\end{Def}

Exact representations of functions by polynomials have been studied
extensively in computer science, where they have numerous applications
in circuit lower bounds\cite{Razborov,Smolensky}, hardness of approximation
\cite{Hastad} and computational learning \cite{MOS}. In these
applications, the set $A$ is generally taken to be $\bp$ or $\zo$.

In this paper, we study a less strict notion of representation of
a function by a polynomial, which is called sign representation. 
\begin{Def}\cite{MP}
\label{def:strong}
A polynomial $P(X_1,\ldots, X_n) \in \R[X_1,\ldots,X_n]$ sign represents $f$ over $A^n$ if
for every $(a_1, \ldots, a_n) \in A^n$, 
\begin{eqnarray*}
f(a_1, \ldots, a_n) = 0 & \Rightarrow & P(a_1, \ldots, a_n) > 0,\\
f(a_1, \ldots, a_n) = 1 & \Rightarrow & P(a_1, \ldots, a_n) < 0.
\end{eqnarray*}
\end{Def}

Such polynomials are also referred to as Polynomial Threshold
Functions for $f$ or Perceptrons. Sign representations have been studied in
computational complexity theory, where they were used by Beigel,
Reingold and Spielman to show that the complexity class \PP\ is
closed under complement \cite{BRS-PP}. Beigel \etal\ use such
representations to show lower bounds on \AC \cite{BRS}. We refer the reader to
the survey by Beigel on applications of such polynomials in complexity
theory\cite{Beigel-survey}. 

Further motivation for studying 
sign representations comes from Valiant's PAC model for computational learning 
 \cite{KV}. If a class of
functions on $n$ variables can be sign represented by degree $d$ polynomials, then that
class can be learnt in time $n^{O(d)}$ in the PAC-learning model (see
\cite{KS} for a precise statement of this result). Indeed, the best
known algorithms for PAC-learning central concept classes like DNF
formulas and intersections of halfspaces use this approach
\cite{KS,KOS}. For this application, there are two parameters of
interest: the degree of the polynomial and the size of its
coefficients. The former determines the running time of the algorithm,
whereas the latter determines the number of samples required
\cite{KS2}. 

While much of the work on computational learning focuses on the Boolean case where
the set $A$ is taken to be $\bp$ or $\zo$, it is quite natural to consider
classes of functions such as decision trees when the variables take
values from larger sets, especially $\{0,\ldots,m\}$. This is explicitly stated as an open problem in
\cite{MOS} for a class of functions called {\em juntas} which are
functions that depend only on some unknown subset of $X_1, \ldots,
X_n$ of size $k$ where $k  = O(\log n)$. It is well-known and easy to
show that juntas are a special class of decision trees \cite{MOS}.

Finally, we will consider an even weaker notion of representing a
function by a polynomial which is called weak representation,
introduced by Aspnes, Beigel, Furst and Rudich.
\begin{Def}
\cite{ABFR}
\label{def:weak}
A polynomial $P(X_1,\ldots, X_n)\in \R[X_1,\ldots,X_n]$ weakly sign represents $f$ over $A^n$ if
for every $(a_1, \ldots, a_n) \in A^n$, 
\begin{eqnarray*}
f(a_1, \ldots, a_n) = 0 & \Rightarrow & P(a_1, \ldots, a_n) \geq 0,\\
f(a_1, \ldots, a_n) = 1 & \Rightarrow & P(a_1, \ldots, a_n) \leq 0.
\end{eqnarray*}
and further $P(X_1,\ldots,X_n)$ does not vanish over the set $A^n$.
\end{Def}
Weak sign representations (or weak representations for short) have been used in computational complexity to
show circuit lower bounds. Aspnes \etal\ \cite{ABFR} and Klivans
\cite{Klivans} use weak representations to show that the parity
function cannot even be approximated by \AC\ circuits. 

Finally, polynomial representations have been studied as a restricted
algebraic model of computation. This gives rise to some natural complexity measures, namely the minimum degree
and sparsity needed to represent a function.
These measures, and tradeoffs
between them have been investigated previously by many researchers
\cite{MP,Beigel-survey,KP,OS-CCC,OS-STOC,KS,KS2,KOS}. Polynomial
representations have also been studied over finite fields
and rings of positive characteristic. This study has yielded useful
insights into computational complexity \cite{Razborov,Smolensky,BBR},
computational learning \cite{MOS} and combinatorics \cite{gro,gro2,CCC}.

\subsection{Our Results}

\begin{Def}
For $A \subset \Z$, the parity function $\Par: A^n \rightarrow \zo$
is defined as 
$$\Par(a_1,\ldots,a_n) = \sum_{i=1}^n a_i \pmod{2}.$$
\end{Def}

We will restrict our attention to the setting where $A$ consists of
non-negative integers, though our methods can be applied to arbitrary
sets $A \subset \Z$. We define the sparsity $\spr(P)$ of a polynomial
$P(X_1,\ldots,X_n)$ to be the number of monomials in its support when
the polynomial is written in the standard monomial basis. 

In this work, we present a systematic study of tradeoffs between
degree and sparsity of sign representations through the lens of the
parity function. 
Our methods also apply to related functions such as
inner-product mod 2 (see Definition \ref{def:ip}). While tradeoffs between degree and sparsity have been investigated
by several researchers \cite{Beigel2,KP,KS2}, previous work focused on
the case $A = \zo$ or $A = \bp$. In contrast, we attempt to prove
bounds that hold for any choice of set $A$. To motivate this, consider
the problem of representing Boolean functions on the $n$-dimensional
hypercube by polynomials. One could identify the hypercube with the set $\{a,b\}^n$
for any $a \neq b \in \R$. Indeed, this freedom to choose the set $A$
is crucially used by the algorithm of Mossel \etal for learning juntas \cite{MOS}.
Thus it is natural to study polynomial representations for  arbitrary
sets $A$. While it is known that the minimum degree of polynomials
representing a function does not depend on the choice of $a$ and $b$,
it is unclear how this affects other parameters
such as coefficient-size and sparsity. 

Obtaining bounds in this general setting is challenging unlike degree, the minimum sparsity
of polynomials representing a function is known to vary greatly with
the choice of the set $A$. However, we show that one can completely
classify the minimum sparsity required to represent parity for any set
$A$ of size $2$. We obtain non-trivial lower bounds on the sparsity
for arbitrary sets $A$ of any size. We obtain tight sparsity lower bounds if we
assume upper bounds on the degree of the polynomial. We show that
there are tradeoffs between the degree and the sparsity of sign
representations. Below we present exact statements of our main results.

We reprove the result of Minsky and Papert \cite{MP} that any polynomial that sign represents parity
over $\zo^n$ has degree $n$ and sparsity $2^n$. We generalize this to show that
representing parity over $\{0,\ldots, m-1\}^n$ with the degree in
each variable at most $m-1$ requires sparsity $m^n$. 
This result shows that low degree representations must have high sparsity.
We show a tradeoff between degree and sparsity by exhibiting sign representations of lower
sparsity but higher degree. We show a lower bound of $n(m-2) + 1$ on the sparsity for
polynomials of any degree representing parity over $\{0, \ldots,
m-1\}^n$. This allows us to prove tight upper and lower bounds for the
case $|A| = 2$. For large sets $A$, we are unable to close the gap
between our upper and lower bounds. 

Our results indicate that studying sparsity gives useful insights into sign representations. 
For instance, consider the polynomials sign-representing $f: \zo^n
\rightarrow \zo$. If we place the restriction that each variable $X_i$ appears with degree at most
$1$, there is a unique polynomial that exactly represents every
function $f$. However even with this restriction, the polynomials that sign
represent a function are not unique: for any $0 < a <b$, the
polynomial $\prod_{i=1}^n(a - bX_i)$
sign represents parity over $\{0,1\}^n$. We show that in any
sign representation of parity, the sign of the coefficient
corresponding to the monomial $\prod_{i \in S}X_i$ must be
$(-1)^{|S|}$, as in the polynomial above. Thus
all sign representations have some similar structure.

As an application of our methods, we show that lower bounds on the sparsity of sign representations
can be used to prove circuit lower bounds for Thresholds of Ands
circuits \cite{Goldmann} (see Definition \ref{def:thr-and}). 
We give a simple proof that any such circuit for parity requires size at least
$(\frac{3}{2})^n$. The best bound known previously was
$(\frac{4}{3})^\frac{n}{2}$. We also show a
lower bound of $2^n$ for computing the inner product function over
$\zo^n \times \zo^n$ which is tight. While our methods are elementary, they give
better lower bounds than those obtained by using the powerful random
restriction method \cite{Goldmann}. Our sparsity-based approach also differs 
from most previous results which related the degree of
sign-representations to the size of small depth
circuits\cite{ABFR,BRS}. Building on our work, Amano and Maruoka
recently used LP-based methods to prove lower bounds on circuits comprising on Thresholds of
symmetric gates that compute the Inner Product function \cite{AM}.

\subsection{Our Techniques}

The degree of sign representations is better understood than the
sparsity for a couple of reasons. Degree is
less dependent than sparsity on the choice of $A$ \cite{MOS}. 
A tool which helps in studying the degree of polynomials sign
representing symmetric functions is symmetrization \cite{MP}: we can assume that the minimum  degree
polynomial sign representing a symmetric Boolean function is symmetric.
Further, if $|A| = m$, we can assume that the minimum degree
polynomial representing a function has degree at most $m-1$ in each
variable. However, such assumptions cannot be made in the context of sparsity.

Our main technical contribution is to show that non-trivial lower
bounds on the sparsity of sign-representations can be obtained using
some elementary techniques and a classical result from algebra called
Descartes' Rule of Signs. Unlike over algebraically closed fields, Descrates' rule of signs shows
that the number of real roots of a univariate real polynomial can be bounded in terms of the number 
of monomials appearing in it (independent of the degree). 
Define the sparsity of a polynomial $P$ to be the number of monomials that occur in it with
non-zero coefficients. We will denote it by $\spr(P)$.

\noindent {\bf Descartes' Rule of Signs: }
{\em Let $P(X)  \in \mathbb{R}[X]$ be a univariate polynomial. Then the
number of positive real roots of $P$ counted with multiplicities is bounded by
the number of sign variations in the sequence of its non-zero coefficients
written in order. In particular, the number of positive roots of $P$
counted with multiplicity is bounded by $\spr(P) - 1$.}

Descartes' rule illustrates that for real univariate polynomials,
sparsity is an important parameter controlling the number of real
zeros. It forms the basis of many
efficient algorithms for real root counting \cite{bpr}. 
An important open problem in real algebraic geometry is to find 
proper analogues of Descartes' rule for multivariate polynomials.
The topological complexity (as measured by the Euler characteristics
or the Betti numbers) of the real zeros of a multivariate real polynomial
can still be bounded in terms of the sparsity of the polynomial independent
of the degree \cite{Khovansky, B99}. However, the known bounds are exponential
in the sparsity and are believed to be nowhere near tight. 
A proper generalization of Descartes' rule to multivariate polynomials
is still elusive and remains a major open problem in real algebraic geometry
(see \cite{Sturmfels} and \cite{Lagarias} for interesting 
conjectures and counter-examples and \cite{Rojas} for results in 
special cases). A small first step in this direction might be to show
tight sparsity bounds for multivariate sign representations of parity for
arbitrary sets $A$.

\subsection{Related Work}

Minsky and Papert prove that representing parity over $\{0,1\}$ inputs needs
degree $n$ and sparsity $2^n$ \cite{MP}. Krause and Pudlak \cite{KP} show that there is a Boolean
function $f$ that has exponential sparsity in the $\{-1,1\}$ basis
but polynomial sparsity in the $\{0,1\}$ basis. O'Donnell and
Servedio \cite{OS-CCC} study various extremal properties of such
representations. The sparsity of random Boolean functions
on $\{-1, +1\}^n$ have been studied in \cite{OS-CCC, Saks}.

%
%
%
%

\section{Preliminaries}

If $P(\Xn) \in \R[X_1, \ldots, X_n]$, we use $P(\Xno, c)$ to denote
the polynomial in $\R[X_1, \ldots, X_{n-1}]$ obtained by substituting $X_n
=c$ in $P(\Xn)$. For $c \in \mathbb{R}$, the sign of $c$ denoted
$\sgn(c)$ is $+1, -1$ or $0$ depending on whether $c$ is positive,
negative or $0$.

The degree of a polynomial $P(\Xn)$ denoted by $\deg(P)$ is the maximum of $\sum
d_i$ over all monomials $ \prod_i X_i^{d_i}$ that occur in the support
of $P(\Xn)$. The degree in the variable $X_i$ which is denoted
$\deg_i(P)$ is  the maximum of $d_i$ over all monomials in the support
of $P(\Xn)$. A multilinear polynomial is one where $\deg_i(P) \leq 1$
for all $i$. The sparsity of a polynomial $P$ denoted $\spr(P)$ is the
number of non-zero monomials in its support.
We also define the sparsity in the variable $X_i$ which we denote
$\spr_i(P)$ to be the number of distinct powers of $X_i$ that occur in $P(\Xn)$. Note
that this is different from the number of monomials in which $X_i$
appears. Given a function $f: A^n \rightarrow \{0,1\}$ define its
complement $\overline{f}:A^n \rightarrow \{0,1\}$ by
$\overline{f}(\an) =  1 - f(\an)$.
If $P(\Xn)$ sign represents $f$, then $-P(\Xn)$ sign
represents $\ol{f}$.

\begin{lemma}
\label{cone}
For $i \in [k]$, let $P_i(\Xn)$ be  polynomials in $\mathbb{R}[\Xn]$
that sign
represent $f$ and let $c_i$ be positive reals. Then
\begin{eqnarray*}
Q(\Xn)= \Sum_{i \in [k]}c_iP_i(\Xn)
\end{eqnarray*}
sign represents $f$.
\end{lemma}
\proof Let $(\an) \in A^n$. Suppose $f(\an) = 0$. Then 
since $P_i(\an) > 0$ for all $i$,  $Q(\an) =
\sum_{i=1}^kc_iP_i(\an) > 0$.
Similarly if $f(\an) =1, Q(\an) < 0$. \qed

Similarly, one can show that if the polynomials $P_i(\Xn)$ weakly sign
represent $f$, then $Q(\Xn)$ also weakly sign represents $f$.

\begin{theorem}
{\bf \cite{polya, bpr} Descartes' Rule of Signs: }
\label{desc} 
Let $P(X) = \sum_{i=0}^nc_iX^i$ be a real univariate polynomial. 
Let $s$ denote the number of sign changes in the sequence $c_0, c_1,
\ldots, c_n$. The number of positive roots of $P(X)$ counted
with multiplicity is bounded by $s$.
\end{theorem}

Let $d_0, \ldots, d_{k-1}$ be non-negative integers such
that $d_0 <  \cdots < d_{k-1}$. Let $a_0, \ldots, a_{k-1}$ be
real numbers such that $a_0 < \cdots < a_{k-1}$. Define the
corresponding generalized Vandermonde matrix as
\begin{equation}
\label{eq:gvd}
V = \left( \begin{array}{llll}
a_0^{d_0} & a_0^{d_1} & \ldots & a_0^{d_{k-1}}\\
a_1^{d_0} & a_1^{d_1} & \ldots & a_1^{d_{k-1}}\\
\ldots & \ldots & \ldots & \ldots\\
a_{k-1}^{d_0} & a_{k-1}^{d_1} & \ldots & a_{k-1}^{d_{k-1}}\\
\end{array}
\right).
\end{equation}
Our goal is to determine the signs of the entries in the inverse of
such a matrix. For this we will use the following lemma:

\begin{lemma}
\cite{polya}
\label{gvd} 
If $a_i > 0$ for all $i$, then $\det(V) > 0$.
\end{lemma}
\proof
The proof is by induction on $k$. The case $k =1$ is trivial.
Assume that the statement holds up to $k - 1$. Now consider the
univariate polynomial in $\R[X]$ defined as
\begin{eqnarray*}
C(X) =  \left| \begin{array}{llll}
a_0^{d_0} & a_0^{d_1} & \ldots & a_0^{d_{k-1}}\\
a_1^{d_0} & a_1^{d_1} & \ldots & a_1^{d_{k-1}}\\
\ldots & \ldots & \ldots & \ldots\\
X^{d_0} & X^{d_1} & \ldots & X^{d_{k-1}}
\end{array}
\right|
\end{eqnarray*}
Let $C(X) =  c_{k-1}X^{d_{k-1}} + c_{k-2}X^{d_{k-2}}+ \cdots + c_0X^{d_0}$.
The sparsity of $C(X)$ is bounded by $k$, hence by Descartes' rule,
it has at most $k-1$ positive roots. But $a_0, \ldots, a_{k-2}$ are
roots of $C(X)$. Hence there are no other roots. Hence the sign at
$a_{k-1}$ (or at any point to the right of $a_{k-2}$) is the same as
the sign at $+\infty$.  This in turn is the sign of the leading coefficient
$c_{k-1}$ of $C(X)$, which is 
\begin{eqnarray*}
c_{k-1} = \left| \begin{array}{lll}
a_{0}^{d_0} & \ldots & a_{k-2}^{d_{k-2}}\\
\ldots & \ldots & \ldots\\
a_{k-2}^{d_0} & \ldots & a_{k-2}^{d_{k-2}}\\
\end{array}
\right|
\end{eqnarray*}
which is positive by the induction hypothesis. \qed

Let $V^{-1} = (v^{-1}_{i,j})$ denote inverse of $V$.
Using Lemma \ref{gvd} and the formula for inverse of a matrix, it is
easy to see that for $0 \leq i, j \leq k -1$, $\sgn(v^{-1}_{ij})= (-1)^{i +j}$.

We will need to consider the case when $a_0 =0$. If $d_0 > 0$, then clearly
the first row is all $0$s and the determinant vanishes. On the other
hand, if $d_0 =0$ we get the matrix
\begin{equation}
\label{eqn:gvd0}
W = 
\left( \begin{array}{llll}
1 & 0 & 0 & 0\\
1 & a_1^{d_1} & \ldots & a_1^{d_{k-1}}\\
1 & \ldots & \ldots & \ldots\\
1 & a_{k-1}^{d_1} & \ldots & a_{k-1}^{d_{k-1}}
\end{array}
\right).
\end{equation}

\begin{lemma}
\label{inv0} 
Let $W^{-1} =(w^{-1}_{i,j})$ denote the inverse of the matrix $W$. For $0 \leq i,j \leq k-1$,
\begin{align*}
\sgn(w^{-1}_{ij})= 
\begin{cases} 
0 & \text{if} \ i = 0, j \geq 1,\\  
(-1)^{i + j} & \text{otherwise.}
\end{cases}
\end{align*}
\end{lemma}
\proof
The minors $W_{ij}$ for $j=0$ and $i \geq 1$ are $0$ since their top
row consists entirely of $0$s. Hence the
entries in $W^{-1}$ for $i =0$ and $j \geq 1$ are $0$. For the other
minors, we can apply Lemma \ref{gvd} to show that they are
positive. Also $\det(W) > 0$, hence by the formula for matrix
inverses,  $\sgn(w^{-1}_{ij}) = (-1)^{i+j}$.  
\qed

%
%
%
%

\section{Lower Bounds}

We first consider the case when $A =\{0,1\}$. Assume that $P(\Xn)$
sign represents parity over $\{0,1\}^n$. If the variable $X$ takes
values in $\{0,1\}$, then $X^k =X$ for $k \geq 2$. So we can use the relation
$X_i^k =X_i$ for $k \geq 2$ to reduce the polynomial $P(\Xn)$ to a
multilinear polynomial. These substitutions can only decrease $\spr(P)$
and $\deg(P)$.


\begin{lemma}
\label{lemma-2}
If $P(\Xn)$ sign represents parity over $\{ 0, 1\}^n$,
\begin{eqnarray}
\label{sum-pq}
P(\Xn) =  X_nQ_1(\Xno) + Q_0(\Xno)
\end{eqnarray}
where $Q_0(\Xno)$ and $-Q_1(\Xno)$ sign represent parity on $n-1$
variables.
\end{lemma}
\proof 
Since $P(\Xn)$ is multilinear, by grouping together monomials
which involve $X_i$, we can write
$$P(\Xn) \ = \ X_nQ_1(\Xno) + Q_0(\Xno)$$
By substituting values for $X_n$, we get  
\begin{align*}
P(\Xno, 0) & = \ Q_0(\Xno), \\
P(\Xno, 1) & = \ Q_1(\Xno) + Q_0(\Xno).
\end{align*}
We now use the so-called self-reducibility of the parity function:
\begin{align*}
\Par(a_1,\ldots,a_{n-1},0) & = \Par(a_1,\ldots,a_{n-1}),\\
\Par(a_1,\ldots,a_{n-1},1) & = \ol{\Par}(a_1,\ldots,a_{n-1}).
\end{align*}
From this it follows that $P(\Xno, 0) = Q_0(\Xno)$ and
$-P(\Xno, 1)$ sign represent parity on $n-1$ variables.
Also, we have 
$$-Q_1(\Xno) \ = \   P(\Xno, 0) - P(\Xno, 1).$$
Hence $Q_1(\Xno)$ sign represents parity by Lemma \ref{cone}. \qed

The polynomial $\prod_{i=1}^{n}(1 - 2X_i)$ sign represents parity over
$\{0,1\}^n$. We will show that the degree and sparsity cannot be lower
for any sign representation.


\begin{theorem}
\cite{MP}
\label{dsp-2}
If $P(\Xn)$ sign represents parity over $\{ 0, 1\}^n$, then
it must have degree $n$ and sparsity $2^n$. 
\end{theorem}
\proof Observe that the sparsity bound of $2^n$ implies that every
monomial including $\prod_{i=1}^nX_i$ has a non-zero coefficient, hence
the degree is $n$. So it is sufficient to prove the sparsity
bound. 

The proof is by induction on $n$.
For $n =1$, let $P(X_1) = aX_1 +b$. $P(X_1)$ must satisfy the conditions
$$P(0)= b > 0, \ \ P(1) = a + b < 0.$$ 
This implies $b > 0$ and $a < -b < 0$, hence $\spr(P) = 2$. 

Assume inductively that the claim holds for $n - 1$ variables.
Write $P(\Xn)$ as in Lemma \ref{lemma-2}. 
Observe that $\spr(P) =\spr(Q_0) + \spr(Q_1)$, since there cannot be cancellations
between the monomials in $X_nQ_1(\Xno)$ and $Q_0(\Xno)$. By
the induction hypothesis $\spr(Q_0) = \spr(Q_1) = 2^{n-1}$, hence $\spr(P)
=2^n$. \qed

We can strengthen the claim to show that the sign of the coefficient
of every monomial is fixed.
For $S \subset [n]$, we denote the coefficient corresponding to the
monomial $\prod_{i \in S}X_i$ by $c_S$. Thus 
$$P(\Xn) = \sum_{S \subseteq [n]}c_S\prod_{i \in S}X_i.$$


\begin{theorem}
\label{sign-2}
If $P(\Xn)$ sign represents parity on $\{0,1\}^n$, then
$\sgn(c_S) = (-1)^{|S|}$.
\end{theorem}
\proof
The proof is by induction on $n$.
The case $n =1$ follows from the Proof of Theorem
\ref{dsp-2}. Assume inductively that the claim holds for $n - 1$
variables. Write $P(\Xn)$ as in Lemma \ref{lemma-2}. The monomials
involving $X_n$ come from $X_nQ_1(\Xno)$ while those not
involving $X_n$ come from $Q_0(\Xno)$. Now consider $S \subset
[n]$ such that $n \notin S$. The coefficient $c_S$ in $P(\Xn)$ is
the same as the coefficient in $Q_0(\Xno)$. Since
$Q_0(\Xno)$ represents parity on $n -1$ variables, hence
$\sgn(c_S) = (-1)^{|S|}$ by the induction hypothesis.  For $S \subset [n]$ such
that $n \in S$, the coefficient $c_S$ in $P(\Xn)$ is equal to the
coefficient $c_{S\setminus \{n\}}$ in $Q_1(\Xno)$. Since $Q_1(\Xno)$
represents the complement of parity, $\sgn(c_S)= -(-1)^{|S| - 1} =
(-1)^{|S|}$ by induction. \qed

One can similarly show a bound on the sum of the coefficient sizes for
polynomials with integer coefficients. We omit the proof.

Next we generalize Theorem \ref{sign-2} to the case when $A = \{0, \ldots, m -1\}$ and the degree in each
variable is at most $m -1$. To construct a polynomial sign representing parity
satisfying these conditions, for $0 \leq j \leq m-2$, let $\alpha_j = j
+ \frac{1}{2}$. Let   
\begin{eqnarray*}
\label{const-m}
P(\Xn) & = & \Prod_{i=1}^n\Prod_{j=0}^{m-2}(-1)^m (X_i - \alpha_j)
\end{eqnarray*}
It can be verified that $P(\Xn)$ indeed sign represents parity on
$A^n$ and $\spr(P) = m^n$. 

Define the univariate polynomial $M(X) \in \R[x]$ by $M(X) = \prod_{j=0}^{m-1}(X -j)$.
Note that $M(X)$ is a monic polynomial of degree $m$ which vanishes on
the set $A = \{0,\ldots, m-1\}$. 
By Euclidean division, for any $d \geq m$, we can write
$$X^d = Q_d(X)M(X) + R_d(X)$$
where $\deg(R_d) \leq m-1$. 

The polynomials $M(X_i)$ for $i \in [n]$ vanish on the set $A^n$.
Given any polynomial $P'(\Xn)$ which sign represents parity over $A^n$,
we can reduce $P'(\Xn)$ modulo the polynomials $M(X_i)$ using 
$$X_i^d \equiv  R_d(X_i)\pmod{M(X_i)}$$
to obtain a polynomial $P(\Xn)$ such that
$$P(\Xn) \equiv P'(\Xn) \pmod{M(X_1),\ldots,M(X_n)}$$
The polynomial $P(\Xn)$ agrees with $P'(\Xn)$ over the set $A^n$, and $\deg_i(P) \leq m-1$. However we will
show that such polynomials where $deg_i(P) \leq m-1$ require sparsity $m^n$.


\begin{lemma}
\label{lemma-m}
Assume that $P(\Xn)$ sign represents parity over  $\{0, 1, \ldots, m
-1\}^n$. If $\deg_n(P) \leq m-1$, then
\begin{eqnarray}
\label{sum-pqm}
P(\Xn) = \Sum_{i=0}^{m-1}X_n^iQ_i(\Xno)
\end{eqnarray}
where $(-1)^{i}Q_i(\Xno)$ represents parity on $n -1$ variables.
\end{lemma}
\proof Since $\deg_n(P) \leq m-1$, grouping monomials by powers of
$X_n$,
\begin{eqnarray*}
P(\Xn) = \Sum_{i=0}^{m-1}X_n^iQ_i(\Xno).
\end{eqnarray*}
Let 
\begin{eqnarray*}
W = 
\left( \begin{array}{llll}
1 & 0 &  \ldots & 0\\
1 & 1 &  \ldots & 1^{m-1}\\
1 & \ldots  & \ldots & \ldots\\
1 & m-1 & \ldots & (m-1)^{m-1}\\
\end{array}\right).
\end{eqnarray*}
By substituting values $0$ through $m -1$ for $X_n$, we get
\begin{eqnarray*}
W\cdot
\left( \begin{array}{l}
Q_0(\Xno)\\
Q_1(\Xno)\\
\ldots\\
Q_{m-1}(\Xno)
\end{array}
\right) & = &  \left( \begin{array}{l}
P(\Xno, 0)\\
P(\Xno, 1)\\
\ldots\\
P(\Xno, m-1)
\end{array}\right) \\
\Rightarrow 
\left( \begin{array}{l}
Q_0(\Xno)\\
Q_1(\Xno)\\
\ldots\\
Q_{m-1}(\Xno)
\end{array}
\right)  & = &  W^{-1}\cdot \left( \begin{array}{l}
P(\Xno, 0)\\
P(\Xno, 1)\\
\ldots\\
P(\Xno, m-1)
\end{array}\right)\nonumber
\end{eqnarray*}
We now expand the LHS. Consider the top row of $W^{-1}$, which is indexed by $i =0$. By Lemma \ref{inv0}, the first
entry is some number $w^{-1}_{00} > 0$, and the other entries are $0$. This implies
$$Q_0(\Xno)  =  w^{-1}_{00}P(\Xno,0),$$
so $Q_0(\Xno)$  sign represents parity on $n-1$ variables. 

For $i \geq 1$, $w^{-1}_{ij} = (-1)^{i +j}|w^{-1}_{ij}|$. Hence
\begin{eqnarray}
Q_i(\Xno) & = & \Sum_{j=0}^{m-1}(-1)^{i +j}|w^{-1}_{ij}|P(\Xno, j)\nonumber\\
\Rightarrow \ \ (-1)^{i}Q_i(\Xno)  & = & \Sum_{j=0}^{m-1}(-1)^{j}|w^{-1}_{ij}|P(\Xno, j). \label{eq:invert}
\end{eqnarray}
We now use the self-reducibility of the parity function:
\begin{align*}
\Par(a_1,\ldots,a_{n-1},j)  = \begin{cases}
\Par(a_1,\ldots,a_{n-1})  & \text{if} \ j \equiv 0 \pmod{2}\\
\ol{\Par}(a_1,\ldots,a_{n-1}) & \text{if} \ j \equiv 1 \pmod{2}
\end{cases}
\end{align*}
Hence the polynomial $(-1)^{j}|w_{ij}|P(\Xno, j)$
sign represents parity on $n -1$ variables for all $j$. Hence by Lemma \ref{cone}, 
$(-1)^iQ_i(\Xn)$ also represents parity on $n-1$ variables for every $i$. \qed


\begin{theorem}
\label{dsp-m}
Let $P(\Xn)$ be a polynomial that sign represents parity over
$\{0, 1, \ldots, m -1\}^n$, with $\deg(X_i) \leq m -1$ for  all $i$. Then
$P(\Xn)$ has sparsity $m^n$ and the sign of the coefficient of the
monomial $\prod_jX_j^{i_j}$ is $(-1)^{\sum_j i_j}$. 
\end{theorem}
\proof 
The proof is by induction.
The base case $n =1$ is an application of Descartes' rule. Let
\begin{eqnarray*}
P(X_1) = \Sum_{i=0}^{m-1}c_iX_1^i
\end{eqnarray*}
Let $0 \leq k \leq m-2$. Since $P(k)$ and $P(k+1)$ have opposite signs, $P(X_1)$ has a root $\alpha_k$
in the interval $(k, k+1)$. Since the degree of $P(X_1)$ is bounded by
$m-1$,
\begin{eqnarray*}
P(X_1) & = & c_{m -1}\Prod_{k=0}^{m-2} (X_1 - \alpha_k)
\end{eqnarray*}
To determine the sign of $c_{m-1}$, substitute $X_1 = 0$.
\begin{eqnarray*}
P(0) & = &  (-1)^{m-1}c_{m-1}\Prod_{k=0}^{m-2} \alpha_k
\end{eqnarray*}
Since $P(X_1)$ represents parity, $P(0) > 0$. Since all the
$\alpha_k$ are positive, we must have $\sgn(c_{m-1}) =
(-1)^{m-1}$. Now applying Descartes' rule, since $P$ has $m -1$
positive roots, there must be $m$ sign changes in the sequence $c_0,
\ldots, c_{m-1}$. Hence, $\sgn(c_i) = (-1)^{i}$. This implies that
$\spr(P) = m$.  

The inductive case proceeds using Lemma \ref{lemma-m} exactly as in
Theorem \ref{dsp-2}. We skip the proof. 
\qed


\begin{corollary}
If $P'(\Xn)$ sign represents parity over  $\{0, 1,
\ldots, m-1\}^n$, then $\deg(P') \geq n(m-1)$.
\end{corollary}
\proof 
We quotient out the $P'(\Xn)$ by the polynomials $M(X_i), \ldots,
M(X_n)$, to get $P(\Xn)$ where $\deg_i(P) \leq m-1$. Note that this
only reduces the total degree, hence $\deg(P') \geq \deg(P)$. By
Theorem \ref{dsp-m}, $\spr(P) \geq m^n$. This implies that for every
tuple $(d_1,\ldots,d_n)$ where $d_i \leq m-1$, the
monomial $\prod_iX_i^{d_i}$ occurs with non-zero
coefficients. Thus the monomial $\prod_{i \in [n]}X_i^{m-1}$ is in the
support, which implies that $\deg(P) \geq n(m-1)$.
\qed

The same proof extends to sets of the form $A =\{a, a+ 1, \ldots, a +
m-1\}$ for $a >0$. This implies the following corollary, by taking $a
=1$ and $m =2$.

\begin{corollary}
\label{par-2}
If $P(\Xn)$ is a multilinear polynomial that sign represents parity over  $\{1, 2\}^n$, it has
sparsity $2^n$.
\end{corollary}

A natural question is what happens to the sparsity if we allow polynomials of higher
degree. It might be that there are polynomials of high degree and low
sparsity and quotienting by the $M(X_i)$s causes the sparsity to
increase. We will address this question in Section $4$. We next turn
our attention to weak representations.

\subsection{Weak Representations}

We first consider weak representations for parity with low
degree. Over $\{0,1\}^n$, the polynomial $P(\Xn) = (-1)^n\prod_iX_i$
gives a weak representation with sparsity $1$, and in fact this is
optimal with regard to degree too. 
\begin{lemma}
\cite{ABFR}
\label{abfr_weak}
Any polynomial that weakly sign represents parity over $\{0,1\}^n$ has
degree $n$. 
\end{lemma}

We show that over $A = \{0, 1, \ldots, m-1\}^n$, a lower bound of $(m-1)^n$
still applies for weak representations when the degree in each
variable is at most $m -1$. 

\begin{lemma}
\label{weak-lemma}
If $P(\Xn)$ weakly sign represents parity over $A = \{0, 1,
\ldots, m-1\}^n$, and if  $\deg_n(P) \leq m-1$, then
\begin{equation}
\label{sum-weak-pqm}
P(\Xn) = \Sum_{i=0}^{m-1}X_n^iQ_i(\Xno)
\end{equation}
where for $i \geq 1$, the polynomial $(-1)^iQ_i(\Xno)$ weakly represents parity on $n -1$
variables.
\end{lemma}
\proof 
The proof is similar to that of Lemma \ref{lemma-m}, the difference being that we need to show
that the polynomials $Q_i(\Xno)$ do not vanish over the set $A^{n-1}$.
By substituting values $0$ through $m -1$ for $X_n$ and inverting the
Vandermonde matrix, we get 
\begin{eqnarray*}
\left( \begin{array}{l}
Q_0(\Xno)\\
Q_1(\Xno)\\
\ldots\\
Q_{m-1}(\Xno)
\end{array}
\right)  =  W^{-1}\cdot \left( \begin{array}{l}
P(\Xno, 0)\\
P(\Xno, 1)\\
\ldots\\
P(\Xno, m-1)
\end{array}\right).\nonumber
\end{eqnarray*}

For $i \geq 1$, by Equation (\ref{eq:invert}), we have $w^{-1}_{ij}
\neq 0$ and
\begin{eqnarray*}
(-1)^{i}Q_i(\Xno)  & = &  \Sum_{j=0}^{m-1}(-1)^{j}|w^{-1}_{ij}|P(\Xno, j).
\end{eqnarray*}

Since $P(\Xn)$ weakly represents parity on $A^n$, the polynomial
$(-1)^jP(\Xno, j)$ either weakly represents parity,
or it vanishes over $A^{n-1}$. Since $P(\Xn)$ is a weak
representation of parity,  it does not vanish on $A^n$. Hence there is
a point $(\an) \in A^n$ so that $P(\an) \neq 0$.
Hence the polynomial $P(\Xno,a_n)$ does not vanish over $A^{n-1}$.
Hence by Lemma \ref{cone},  the polynomials $(-1)^iQ_i(\Xno)$ weakly sign represent
parity on $A^{n-1}$ for $i \geq 1$. \qed

The condition $i \geq 1$ in the statement of Lemma
\ref{weak-lemma} is in fact necessary: take the polynomial $P(\Xn) =
(-1)^n\prod_iX_i$ that weakly sign represents parity on $\zo^n$.
In this case, $Q_0(\Xno) = 0$, so it does not represent parity even weakly.

We use Lemma this to show a lower bound of $(m-1)^n$ on the sparsity of
weak representations over $A^n$. The base case $n =1$ is proved using
Lemma \ref{abfr_weak}.  

\begin{lemma} 
\label{our_weak}
Any univariate polynomial $P(X) \in \R[X]$ that weakly sign represents parity
over $\{0, \ldots, m-1\}$ must have $m - 2$ roots in the interval $(0,
m -1]$. 
\end{lemma}
\proof 
We first show that $\deg(P) \geq m-1$. Assume that this is not so.
Let $X_1, \ldots, X_{m-1}$ be variables that take values in $\{0,1\}$. Then the polynomial
$$Q(X_1,\ldots,X_{m-1})= P(\Sum_{i=1}^{m-1} X_i)$$
weakly sign represents parity on
$\{0,1\}^{m-1}$ and $\deg(Q) < m-1$, contradicting
Lemma \ref{abfr_weak}. 

Now consider the factorization of $P(X)$ over the reals. Assume that
this contains an irreducible polynomial $D(X) \in R[X]$ with $\deg(D) =2$. 
The polynomial $D(X)$ does not have real roots 
its sign stays unchanged in the interval $[0,\ldots,m-1]$. Hence we can replace $D(X)$
by the constant $\sgn(D(0))$, and get a sign representation of lower degree. 
Similarly, consider a linear factor of the form $(X -
\alpha)$ where $\alpha \notin [0, m-1]$. Such linear factors can also
be replaced by their signs at $0$. Further we may assume that there is a root of multiplicity at most $1$
at $0$. If not, we can write $P(X) = X^kQ(X)$ for $k \geq 2$. The polynomial
$XQ(X)$ has the same sign at each point in $[0, m-1]$ and only smaller degree. 
We are left with a polynomial of the form 
$$P'(X) = \prod (X - \alpha_i) \hspace{1cm} \alpha_i \in [0, m-1] $$
which weakly represents parity over $\{0, \ldots, m-1\}$, hence $\deg(P') \geq m-1$.
Since $0$ is a root of multiplicity at most $1$, at least $m-2$ of the
roots $\alpha_i$ lie in the interval $(0, m-1]$. \qed

\begin{theorem}
\label{weak-lower}
Let $P(\Xn)$ be a polynomial that weakly represents parity over
$\{0, \ldots, m-1\}^n$. If $\deg_i(P) \leq m-1$ for every $i \in [n]$, then $\spr(P) \geq m^n$.
\end{theorem}
\proof The proof is by induction on $n$. When $n =1$, by Lemma \ref{our_weak}
$P(X_1)$ has $m-2$ roots in $(0,m-1]$, hence by Descartes' rule,
$\spr(P) \geq m-1$. 

For the inductive case, we use Lemma \ref{weak-lemma}.
From Equation (\ref{sum-weak-pqm}) it follows that
$$\spr(P) \geq \sum_{i=1}^{m-1} \spr(Q_i).$$
By Lemma \ref{weak-lemma}, the polynomial $(-1)^iQ_i(\Xno)$ weakly represents parity on $A^{n-1}$.
Hence by induction, $\spr(Q_i) \geq (m-1)^{n-1}$, and so $\spr(P) \geq (m-1)^n$. \qed

This bound is in fact tight.

\begin{lemma}
\label{lem:example}
There is a polynomial $P(\Xn)$ that weakly represents parity over
$\{0, \ldots, m-1\}^n$ where $\deg_i(P) \leq m-1$ for every $i \in [n]$,
and $\spr(P) = (m-1)^n$.
\end{lemma}
\proof
Take $Q(\Xn)$ to be a polynomial that sign represents parity on $\{1,\ldots,m-1\}^n$ satisfying $\spr(Q) =
(m-1)^n$, and $\deg_i(Q) \leq m-2$. We claim that the polynomial 
$$P(\Xn) = Q(\Xn)\cdot \Prod_{i=1}^nX_i$$ 
weakly represents parity over $\{0,\ldots, m-1\}^n$. This is
because, for $\an \in \{0,\ldots, m-1\}^n$,
$$\sgn(P(\an)) = \begin{cases} 
0 & \text{if } a_i = 0 \text{ for some $i$},\\
\sgn(Q(\an)) & \text{otherwise.}
\end{cases}$$
Further, since $\deg_i(Q) \leq m-2$, $\deg_i(P) \leq m-1$ for all $i
\in [n]$. Also, $\spr(P) = \spr(Q) = (m-1)^n$.
\qed

The proof of Lemma \ref{lem:example} crucially uses the fact that $0
\in A$. Indeed we will show that if $A =\{1, \ldots, m\}$, then
weak representations of parity require sparsity $m^n$.

\begin{corollary}
\label{cor:weak}
Let $P(\Xn)$ be a polynomial that weakly represents parity over
$\{1, \ldots, m\}^n$. If $\deg_i(P) \leq m-1$ for every $i \in [n]$, then $\spr(P) \geq m^n$.
\end{corollary}
\proof 
Let $P(\Xn)$ be as above. The polynomial 
$$P'(\Xn) = P(\Xn)\cdot \Prod_{i=1}^nX_i$$
weakly represents parity over $\{0,\ldots, m\}^n$. 
Further, since $\deg_i(P) \leq m-1$, $\deg_i(P') \leq m$ for all $i \in [n]$.
Hence we can apply Theorem \ref{weak-lower}, which implies
$$\spr(P) = \spr(P') \geq m^n.$$ 
\qed

%
%
%

\section{Upper Bounds}
Does the lower bound of $m^n$ in Theorem \ref{dsp-m} hold for all
polynomials? Or are there polynomials with higher degree but lower sparsity?
We show that such a tradeoff is indeed possible.

\begin{theorem}
\label{upper-12} 
There exists a polynomial $P(\Xn)$ that sign
represents parity over $\{1, 2\}^n$ with $\deg(P) = n^2$ and $\spr(P) = n + 1$.
\end{theorem}
\proof
Define $w: \{1,2\}^n \rightarrow \{2^k\}_{k=0}^n$ by $w(\an) = \prod_ia_i$. 
If $w(\an) = 2^k$, then $a_i = 2$ for exactly $k$
co-ordinates $i \in [n]$ hence $\Par(\an) \equiv n -k \pmod{2}$.

Choose points $\alpha_j \in (2^{j-1}, 2^j)$ for $1 \leq j \leq n$ and let
$$P(\Xn) = \Prod_{j=1}^n\left(\prod_{i=1}^nX_i - \alpha_j\right),$$
so that $\deg(P) = n^2$ and $\spr(P) = n +1$. We claim that
$P(\Xn)$ sign represents parity on $\{1,2\}^n$.
Note that
$$P(\an) \ = \Prod_{j=1}^n(w(\an) - \alpha_j).$$
If $w(\an) = 2^k$, then
\begin{align*} 
P(\an) & = \Prod_{j=1}^n( 2^k - \alpha_j)\\
\Rightarrow \ \sgn(P(\an)) & = (-1)^{n-k}.
\end{align*}
Thus the polynomial $P(\Xn)$ sign represents parity. Its
sparsity is $n +1$ and its degree is $n^2$. \qed

In contrast, Corollary \ref{par-2} shows a lower bound of $2^n$ on the
sparsity for sign representations by multilinear polynomials.

We can extend Theorem \ref{upper-12} to show that for any set $A$ of
non-negative integers of size $m$, there are polynomials that weakly
sign represent parity whose sparsity is less $m^n$, but which have
high degree. 

\begin{theorem}
\label{weak-upper}
For any set $A$ of non-negative integers of cardinality $m$, parity
can be weakly sign represented over $A^n$ by a polynomial $P(\Xn)$ that has sparsity at most
$\binom{n +m -1}{n}$ 
and degree $n\cdot \binom{n + m -1}{n}$.
\end{theorem}
\proof Define the function $w: A^n \rightarrow \Z$ by $w(\an) =
\prod_i a_i$. This maps $A^n$ to a set $S$ of size at most
$\binom{n +m -1}{n}$ in $\mathbb{Z}$. Let $a$ denote the largest
integer in the set $A$. Note that $w(a,\ldots,a) = a^n$ is the largest
integer in $S$. Further
$(a,\ldots,a)$ is the unique point in $A^n$ that is mapped to $a^n$ by
$w$. We claim that the polynomial 
$$P(\Xn) = (-1)^{na}\cdot \prod_{\alpha \in S \setminus\{ a^n\}}\left( \prod_{i \in
    [n]}X_i - \alpha \right)$$ 
weakly represents parity on $A^n$. To prove this, note that
$$P(\an) = (-1)^{na}\cdot \prod_{\alpha \in S \setminus\{ a^n\}}(w(\an) - \alpha).$$ 
Hence the polynomial vanishes for every point in $A^n \setminus \{
(a,\ldots,a)\}$. Further
\begin{align*}
P(a,\ldots,a) & = (-1)^{na}\cdot \prod_{\alpha \in S \setminus\{ a^n\}}(a^n - \alpha)\\
\Rightarrow \sgn(P(a,\ldots,a)) & = (-1)^{na}.
\end{align*}
Thus $P(\Xn)$ weakly represents parity on $A^n$.
\qed

In general $\spr(P)= |S|$ and $|S|$ depends on the set $A$. 
For some sets $A$, $\spr(P)$ can be significantly smaller that the bound stated $\binom{n
  +m -1}{n}$ stated above.  
In the case when $A = \{0,\ldots, m-1\}$ and $n$ is a fixed constant, one can in fact show that
$\spr(P) = o(m)^n$. This is a
consequence of Erd\"{o}s' multiplication table theorem which states that
the number of distinct integers less than $m^n$ which can be expressed
as the product of $n$ numbers each less than $m$ is $o(m)^n$ 
\cite{erdos}.

\section{Lower Bounds without Degree Restrictions}
We will now show a lower bound which holds for all polynomials
strongly representing parity on $\{1, \ldots, m\}^n$ without any
restrictions on the degree or sparsity of each variable. The proof is
a generalization of the proof idea of Theorem \ref{dsp-m}.

\begin{theorem}
\label{general-lower} Let $P(\Xn)$ be a polynomial
which sign represents parity over $\{1, \ldots, m\}^n$. Then $\spr(P)
\geq n(m -1) + 1$.
\end{theorem}
\proof The proof is by induction on $n$. When $n = 1$, the claim follows
by Descartes' rule. Assume it is true for $n -1$. Recall that
$\spr_n(P)$ is the number of distinct powers of $X_n$ that occur in
monomials in the support of $P$. If we set all the other variables to $1$, the
univariate polynomial $Q(X_n) = (-1)^{n-1}P(1,\ldots,1,X_n)$ sign represents parity on $\{1,
\ldots, m\}$, hence it must have sparsity at least $m$. Hence $\spr_n(P)
\geq \spr(Q) \geq m$. If $k > n(m-1) + 1$, there is nothing to prove. Hence we may assume $m \leq k
\leq n(m-1)$. 

Grouping monomials in $P(\Xn)$ by the power of $X_n$ they contain, we
can write
\begin{eqnarray*}
P(\Xn) & = & \Sum_{i=1}^kX_n^{d_i}Q_i(\Xno).
\end{eqnarray*}
By substituting values $1$ through $m$ for $X_n$, we get
\begin{equation}
\label{eq:matrix}
\left( \begin{array}{llllll}
1 & \cdot \cdot & 1 & 1 & \cdot \cdot & 1\\
2^{d_1} & \cdot \cdot & 2^{d_m} & 2^{d_{m +1}} & \cdot \cdot & 2^{d_k}\\
\cdot \cdot & \cdot \cdot & \cdot \cdot  & \cdot \cdot & \cdot \cdot
& \cdot \cdot\\
m^{d_1} & \cdot \cdot & m^{d_m} & m^{d_{m+1}} & \cdot \cdot & m^{d_k}\\
\end{array}\right)
\left( \begin{array}{l}
Q_1(\Xno)\\
Q_2(\Xno)\\
\ldots\\
Q_k(\Xno)
\end{array}
\right) \
= \  \left( \begin{array}{l}
P(\Xno, 1)\\
P(\Xno, 2)\\
\ldots\\
P(\Xno, m)
\end{array}\right)
\end{equation}
We denote the $m \times k$ matrix by $A$. While we cannot prove
that each $Q_i(\Xno)$ represents parity (or its complement), we will show that
appropriate linear combinations of the $Q_i(\Xno)$ sign represent
parity. We pre-multiply each side of Equation (\ref{eq:matrix}) by $U$, which is the inverse
of the $m \times m$ generalized Vandermonde matrix consisting of
the first $m$ columns of $A$.
\begin{eqnarray*}
\left( \begin{array}{llllll} 1 &  \cdot \cdot & 0 & b_{1, m+1}
& \cdot \cdot & b_{1, k}\\
0 & \cdot \cdot & 0 & b_{2, m+1} & \cdot \cdot& b_{2, k}\\
\cdot \cdot & \cdot \cdot & \cdot \cdot & \cdot \cdot & \cdot \cdot&
\cdot \cdot \\
0 &  \cdot \cdot & 1 & b_{m, m+1} & \cdot \cdot & b_{m, k}\\
\end{array} \right)
\left( \begin{array}{l}
Q_1(\Xno)\\
Q_2(\Xno)\\
\ldots\\
Q_{k}(\Xno)
\end{array}
\right) &
=  & U\cdot \left( \begin{array}{l}
P(\Xno, 1)\\
P(\Xno, 2)\\
\ldots\\
P(\Xno, m)
\end{array}\right)
\end{eqnarray*}
Using the sign alternations of the entries of $U$, we
conclude that for $1 \leq i \leq m$ the polynomials
\begin{eqnarray*}
R_i(\Xno) & = & Q_i(\Xno) + \Sum_{j=m+1}^kb_{ij}Q_j(\Xno)
\end{eqnarray*}
sign represent parity or its complement on $n-1$ variables. 
Hence by applying the induction hypothesis, 
$$\spr(R_i) \geq (n-1)(m-1) +1.$$ 
But we also have
$$\spr(R_i) \leq \spr(Q_i) + \sum_{j=m+1}^k \spr(Q_j).$$ 
Hence we get
$$\spr(Q_i) + \sum_{j=m+1}^k \spr(Q_j) \ \geq \ (n-1)(m-1) + 1.$$
By choosing the matrix $U$ to be the inverse of an appropriate
sub-matrix, we can obtain a similar equation for any subset of the
$Q_i$s of cardinality $k - m + 1$. There are $\binom{k}{m -1}$ such
subsets. Each $Q_i$ occurs in exactly $\binom{k -1}{m - 1}$ of them. Hence we get
\begin{eqnarray*}
\binom{k -1}{m -1}\Sum_{i=1}^k \spr(Q_i)  & \geq & \binom{k}{m -1}((n-1)(m-1) + 1)\\
\Rightarrow \Sum_{i=1}^k \spr(Q_i) & \geq  & \frac{k}{k -m +1}((n-1)(m-1) + 1)
\end{eqnarray*}
The quantity $\frac{k}{k - m +1}$ monotonically decreases as $k$
increases. In the range $m \leq k \leq n(m -1)$, it is always
greater than $\frac{n(m-1) +1}{(n-1)(m-1) + 1}$ which is the value
it takes for $k = n(m-1) + 1$. Hence
\begin{eqnarray*}
\sum_{i=1}^k \spr(Q_i) & \geq & n(m-1) + 1
\end{eqnarray*}
But $\spr(P) = \sum_{i=1}^k \spr(Q_i)$, hence the claim is proved. \qed

\begin{corollary}
Any polynomial that sign represents parity over $\{1, 2\}^n$ must
have sparsity at least $n + 1$.
\end{corollary}
This follows by substituting $m =2$ in Theorem
\ref{general-lower}. This shows that the construction of
Theorem \ref{upper-12} is optimal with regard to sparsity. 
We can now prove tight lower bounds on polynomials sign
representing parity on $A^n$ for any set $A$ of size $2$.
Let $A =\{a,b\}$ where $0 \leq a < b$.
\begin{itemize}
\item If $a =0$, then any polynomial which sign represents parity has
  sparsity at least $2^n$.
\item If $a > 0$, then any polynomial which sign represents parity has
  sparsity at least $n + 1$.
\end{itemize}

While the lower bound of $n(m-1) + 1$ in Theorem
\ref{general-lower} is tight for $m =2$, this is far from the
upper bound of Theorem \ref{upper-12} for large $m$. It would be
interesting to close this gap.

%
%
%
%
\section{Circuit Lower Bounds}
We shall use bounds on the sparsity of parity to derive lower bounds on the
size of certain restricted circuits. The circuits we consider are
rather weak, however the proof of the lower bound is simple and yields
better parameters than were previously known. 

\begin{definition} 
A Threshold function $f:\zo^n \rightarrow \{0,1\}$ is defined as
\begin{align*}
f(\an) = \begin{cases} 1 & \text{if} \ w_0 +\sum_{i=1}^nw_ia_i < 0\\
0 & \text{if} \ w_0 +\sum_{i=1}^nw_ia_i > 0\\
\end{cases}
\end{align*}
where $w_0,\ldots, w_n \in \R$. The coefficients $w_i$ are called the
weights of the Threshold function. A gate computing a threshold
function is called a Threshold gate and is denoted by \Th.
\end{definition}
In our definition, we assume that
$w_0,\ldots,w_n$ are such that $w_0 +\sum_{i=1}^nw_ia_i \neq 0$ for $\an \in \zo^n$.

\begin{definition}
\label{def:thr-and}
A Threshold of Ands circuit (denoted \Thr ) is a depth-two circuit with Boolean inputs $X_1,
\ldots, X_n$. The top level of the circuit consists of a single
\Th\ gate, while the bottom level consists of \AND\ gates. The
inputs to the \AND\ gates are the inputs $X_1, \ldots, X_n$ and their
complements $\neg X_1, \ldots, \neg X_n$. The size of a \Thr\ circuit is
defined to be the number of \AND\ gates in the bottom level of the
circuit. The minimum size of a \Thr\ circuit needed to compute a Boolean
function $f$ is denoted by $S(f)$.
\end{definition}

These circuits are
well-studied (see \cite{Goldmann} and the references therein). By
De Morgan's law, such circuits can simulate \Or\ gates at the bottom level.
We will show that $S(f)$ corresponds to minimum sparsity required to
sign represent $f$ over a certain {\em  basis}. Thus proving circuit
lower bounds is equivalent to proving bounds on the sparsity of sign
representations. 

To begin with, assume that the inputs to the \AND\ gates were only
the variables $X_i$, not their complements.
Each \AND\ gate computes a function of the form $\prod_{i \in A} X_i$ where $A$ is the set of inputs into the
gate. Such a circuit computing parity corresponds to a sign
representation of parity in the standard monomial basis. The number of \AND\
gates is exactly the number of non-constant monomials required. By
Theorem \ref{dsp-2}, this is $2^n  -1$. 

In a general \Thr\ circuit, and $\AND$ gate computes the function
$\wedge_{i \in I} X_i \wedge_{j \in J} \neg X_j$. We can assume that $I
\cap J$ is empty, else the \AND\ gate computes the function $0$.
Thus the  \AND\ gate computes the polynomial
\begin{equation}
\label{and-gate}
B(\Xn) = \prod_{i \in I}X_i \prod_{j \in J}(1 - X_j), \ \ I \cap J = \phi 
\end{equation}
Let $\mc{B}_n$ denote the set of all such polynomials taken over all
choices of the sets $I$ and $J$. It is easy to
show that $|\mc{B}_n| = 3^n$. Since $\mc{B}_n$ contains the standard
monomial basis, so it spans the $\R$-vector space of multilinear
polynomials in $\R[X_1,\ldots,X_n]$. Since this vector space has
dimension $2^n$, there are many ways to write a
multilinear polynomial $P(\Xn)$ as a linear combination of polynomials
in $\mc{B}_n$. We will define $\sprb(P)$ as the minimum
possible sparsity over all such linear combinations. Formally:

\begin{definition}
Let $P(\Xn)$ be a multilinear polynomial in $\R[\Xn]$. We define the
sparsity of $P(\Xn)$ over $\mc{B}_n$ as
$$\sprb(P) = \{\min k \ |\ P(\Xn) = \sum_{i=1}^k c_iB_i(\Xn), \ \ B_i(\Xn)
\in \mc{B}_n\}.$$ 
\end{definition}

The following lemma relating circuit-size for \Thr\ circuits computing
$f$ and sparsity over $\mc{B}$ of polynomials that sign represent $f$
follows from the preceding discussion.

\begin{lemma}
\label{size-sparsity}
For any Boolean function $f:\{0,1\}^n \rightarrow \{0,1\}$,
$$S(f) = \min \sprb(P)$$
over all polynomials $P(\Xn)$ that sign represent $f$.
\end{lemma}

\begin{theorem}
\label{size-parity}
Every \Thr\ circuit computing the parity function on $\zo^n$ has 
size at least $\left( \frac{3}{2} \right)^n$.
\end{theorem}
\proof 
We will show that if $P(\Xn)$ sign represents parity, then $\sprb(P)
\geq (\frac{3}{2})^n$.

The proof is by induction on $n$. For $n =1$, $\mc{B}_1 = \{X_1, 1
-X_1, 1\}$. Since none of these polynomials or their multiples sign
represents parity on $1$ variable, $\sprb(P) \geq 2$.

Now assume the claim holds for
$n -1$. Let $P(\Xn)$ sign represent parity.  Consider the
sparsest representation of $P$ over $\mc{B}$.
\begin{eqnarray*}
P(\Xn) & = & \Sum_{I\cap J = \phi}c_{I,J}\Prod_{i \in I}X_i\Prod_{j \in
  J}(1 -X_j)
\end{eqnarray*}
Grouping together monomials where $X_n$ appears, monomials where
$(1 -X_n)$ appears, and those where neither appears, we get
\begin{equation}
\label{eq:abc}
P(\Xn) \ = \   X_nA(\Xno) + (1 - X_n)B(\Xno) + C(\Xno).
\end{equation}
The best (sparsest) way to write $P(\Xn)$ as a linear combination of
polynomials in $\mc{B}_n$ is to use the best (sparsest) expression for each of $A(\Xno),
B(\Xno)$ and $C(\Xno)$ as linear combinations of polynomials
in $\mc{B}_{n-1}$. Hence 
$$\sprb(P) = \sprb(A) + \sprb(B) + \sprb(C).$$
Substituting for $X_n$ in Equation (\ref{eq:abc}), 
\begin{eqnarray*}
P(\Xno, 0)  & =   B(\Xno) + C(\Xno),\\ 
P(\Xno, 1)  & =   A(\Xno) + C(\Xno),\\
P(\Xno, 0)  -   P(\Xno, 1) & =  B(\Xno) - A(\Xno)
\end{eqnarray*}
All the polynomials on the LHS represent either parity or its
complement on $n -1$ variables. By applying the induction hypothesis, 
\begin{eqnarray*}
\sprb(B) + \sprb(A) \ \geq \ \sprb(B - A) \ > \ (3/2)^{n-1},\\
\sprb(B) + \sprb(C) \ \geq \ \sprb(B + C) \ > \ (3/2)^{n-1},\\
\sprb(A) + \sprb(C) \ \geq \ \sprb(A + C) \ > \ (3/2)^{n-1}.
\end{eqnarray*}
Adding these equations, we get
\begin{eqnarray*}
\sprb(P) = \sprb(A) + \sprb(B) + \sprb(C) > (3/2)^n
\end{eqnarray*}  
which completes the proof of the Theorem.
\qed

\begin{proposition}
\label{beigel}
There is a \Thr\ circuit computing the parity function on $\zo^n$ of size $O(5^{\frac{n}{3}})$.
\end{proposition}
\proof The polynomial 
\begin{eqnarray*}
Q(X_1,X_2,X_3) & = & X_1X_2X_3 + X_1(1-X_2)(1-X_3) + 
X_2(1-X_3)(1-X_1)+X_3(1-X_1)(1-X_2)
\end{eqnarray*}
exactly represents parity on $\zo^3$.
Hence, the polynomial
$P(X_1,X_2,X_3) = 1 - 2Q(X_1,X_2,X_3)$ sign represents parity
on $3$ variables and has sparsity $5$.
Hence the polynomial 
$$R(X) = \prod_{i=1}^{\frac{n}{3}}P(X_{3i -2}, X_{3i -1}, X_{3i})$$
sign represents parity on $n$ variables, and $\sprb(R) = 5^{\frac{n}{3}}$. 
\qed

\begin{Def}
\label{def:ip}
The inner product function $\IP: \zo^n\times \zo^n \rightarrow \zo$
is defined as 
$$\IP(a_1,\ldots,a_n,b_1,\ldots,b_n) = \sum_i a_ib_i \pmod{2}.$$
\end{Def}

\begin{theorem}
\label{size-ip}
Every \Thr\ circuit computing the inner product function on
$\zo^n\times \zo^n$ has size at least $2^n$.
\end{theorem}
\proof The proof is by induction on $n$. The base case is trivial.
Assume the claim holds for $n-1$. 

Let $P(\Xn,\Yn)$ sign represent \IP\ on $\zo^n\times \zo^n$.
Consider the sparsest way to write $P(\Xn,\Yn)$ over $\mc{B}_{2n}$. Grouping the monomials according to
$X_n,Y_n$, where each $A_i$ is a polynomial in $X_1, Y_1, \ldots,
X_{n-1}, Y_{n -1}$, 
\begin{eqnarray*}
P(\Xn, \Yn) \ =  X_nY_nA_1 + X_n(1-Y_n)A_2  + (1 -
X_n)Y_nA_3 + (1 - X_n)(1 - Y_n)A_4 \\
+ X_nA_5 + Y_nA_6 + (1 -X_n)A_7 + (1 - Y_n)A_8 + A_9
\end{eqnarray*}
Now substituting for $X_n, Y_n$ and writing $P(0,1)$ for $P(\Xno, 0,\Yno, 1)$ and so on,
\begin{eqnarray}
P(0, 1) &  =   & A_3 + A_6 + A_7 + A_9\label{eq:01},\\
P(1, 0) &  =   & A_2 + A_5 + A_8 + A_9\label{eq:10},\\
P(1, 1) &  =   & A_1 + A_5 + A_6 + A_9\label{eq:11}.
\end{eqnarray}
Subtracting Equation (\ref{eq:11})  from (\ref{eq:01}) and
(\ref{eq:10}) respectively, 
\begin{eqnarray}
P(0, 1) - P(1,1) &  =   & - A_1 + A_3 - A_5 + A_7, \label{eq:sub1}\\
P(1, 0) - P(1,1) &  =   & - A_1 + A_2 - A_6 + A_8. \label{eq:sub2}
\end{eqnarray}
In the above equations the polynomials on the LHS represents \IP\ or its complement on
$\zo^{n-1}\times\zo^{n-1}$, so each has sparsity at least $2^{n-1}$ by the
induction hypothesis. Applying this observation to equations
\ref{eq:01}, \ref{eq:10}, \ref{eq:sub1} and \ref{eq:sub2},
\begin{eqnarray*}
\sprb(A_3) + \sprb(A_6) + \sprb(A_7) + \sprb(A_9)  \geq 2^{n-1},\\
\sprb(A_2) + \sprb(A_5) + \sprb(A_8) + \sprb(A_9)  \geq 2^{n-1},\\
\sprb(A_1) + \sprb(A_3) + \sprb(A_5) + \sprb(A_7)  \geq 2^{n-1},\\ 
\sprb(A_1) + \sprb(A_2) + \sprb(A_6) + \sprb(A_8)  \geq 2^{n-1}.
\end{eqnarray*}
Adding these equations, we get
$$2\left(\Sum_{i =1}^9 \sprb(A_i)\right) - 2\sprb(A_4)\geq 4\cdot 2^{n-1}.$$
Hence 
$$\sprb(P) = \Sum_{i=1}^9\sprb(A_i) \geq 2^n$$
which completes the proof.
\qed

\begin{proposition}
There is a \Thr\ circuit computing the inner product function on $\zo^n\times
\zo^n$ of size $2^n$.
\end{proposition}
\proof 
The polynomial 
\begin{eqnarray*}
Q(X_1,X_2, Y_1,Y_2) = X_1Y_1 + X_2Y_2 - 2X_1Y_1X_2Y_2
\end{eqnarray*}
exactly represents $\IP$ on $\zo^2\times \zo^2$. The polynomial
$$P(X_1,X_2, Y_1,Y_2)= 1 - 2Q(X_1,X_2 Y_1,Y_2)$$ 
sign represents $\IP$ and on $\zo^2\times \zo^2$ and  $\sprb(P) =
4$. For arbitrary $n$, the polynomial
$$R(\Xn,\Yn)=  \prod_{i=1}^{\frac{n}{2}}P(X_{2i-1},X_{2i},Y_{2i -1},Y_{2i})$$
sign represents \IP\ on $\zo^n\times \zo^n$  and $\sprb(R) = 4^{\frac{n}{2}} = 2^n$.
\qed

\section*{Acknowledgments}
We would like to thank Richard Beigel for pointing out the upper bound
of Proposition \ref{beigel}, as well as many pointers to
literature. We would like to thank Ernie Croot for the reference to
Erd\"{o}s' multiplication table theorem \cite{erdos}. We thank the
anonymous referee for numerous suggestions that helped improve the
presentation.

\bibliographystyle{alpha}
\bibliography{sparse2}

\end{document}